\documentclass[11pt]{amsart}
\usepackage{amssymb, amsthm, amsmath}
\usepackage{mathtools}
\usepackage{mathrsfs}
\usepackage{hyperref, manfnt}
\usepackage[numbers]{natbib}
\usepackage{setspace}
\usepackage{graphicx}
\usepackage{enumitem}
\usepackage{MnSymbol}
\usepackage{fullpage}
\usepackage{tikz}
\usepackage{tikz-cd}
\newcommand*{\DashedArrow}[1][]{\mathbin{\tikz [baseline=-0.25ex,-latex, dashed,#1] \draw [#1] (0pt,0.5ex) -- (1.3em,0.5ex);}}%
\usepackage{wrapfig}
\usepackage{float}
\usepackage{microtype}
\makeatletter
\newcommand*\bigcdot{\mathpalette\bigcdot@{.5}}
\newcommand*\bigcdot@[2]{\mathbin{\vcenter{\hbox{\scalebox{#2}{$\m@th#1\bullet$}}}}}

\makeatother
\newcommand{\RomanNumeralCaps}[1]
{\MakeUppercase{\romannumeral #1}}
\usepackage{kantlipsum}
\newcommand\Z{\mathbb{Z}}
\usepackage{bm}
\usepackage[cmtip,all]{xy}

\newtheorem*{theorem*}{Theorem}

\newtheorem*{defi*}{Definition}
\numberwithin{theorem}{section}

\newtheorem{rem}{Remark}
\newtheorem*{cor*}{Corollary}
\newtheorem{lemma}{Lemma}
\newtheorem*{ack}{Acknowledgements}

\newtheorem*{rem*}{Remark}
\numberwithin{lemma}{section}
\newtheorem*{lemma*}{Lemma}

\newtheorem*{claim*}{Claim}

\tikzset{
	symbol/.style={
		draw=none,
		every to/.append style={
			edge node={node [sloped, allow upside down, auto=false]{$#1$}}}
	}
}

\theoremstyle{definition}

\newcommand{\p}{\mathbb{P}}

\usepackage{stackengine}
\newcommand{\hollowslash}{\setbox0=\hbox{/}\def\holwd{3pt}%
  \stackengine{-.3pt}{/}{\rlap{\kern.5pt\rule{\holwd}{0.4pt}}}{O}{r}{F}{F}{S}%
  \kern\dimexpr\holwd-\wd0-.2pt\relax%
  \stackengine{-.4pt}{/}{\llap{\rule{\holwd}{0.4pt}\kern-0.1pt}}{U}{l}{F}{F}{S}%
}
\begin{document}
\renewcommand*{\thefootnote}{\alph{footnote}}
\title[Unirationality of the universal moduli of semistable bundles over curves]{\textbf{Unirationality of the universal moduli space of\\ semistable bundles over smooth curves}}
\author{Shubham Saha}
\address{Department of Mathematics, University of California, San Diego, La Jolla, CA 92093, USA}
\email{shsaha@ucsd.edu}
\begin{abstract}\noindent
We construct explicit dominant, rational morphisms from projective bundles over rational varieties to relevant moduli spaces, showing their unirationality. These constructions work for $U(r,d,g)$; for all ranks, degrees and genus $2\leq g \leq 9$. Furthermore, the arguments presented also show that a similar conclusion can be made for $U(r, \mathcal{L}, g)$ for all $r,d$ and unirational $M_g$. 
\end{abstract}
\maketitle{}
\section{Introduction}
The geometry of a moduli space is greatly influenced by its unirationality. The moduli space of curves of genus $g$ over the complex field $\mathbb{C}$ has been the subject of multiple investigations throughout the years and its unirationality has been answered almost entirely in \cite{ac}, \cite{cr1}, \cite{ser}, \cite{sev}, \cite{verra}, \cite{EH} and \cite{FJP}, other than values $15\leq g\leq 21$. The moduli space of semistable bundles with fixed rank and determinant over a given smooth projective curve has also been the subject of thorough investigation in regards to their rationality in \cite{N} and \cite{SK}. We shall be considering the universal moduli space of semistable bundles over smooth curves and universal moduli space of semistable bundles over smooth curves with fixed determinant, first compactified by R. Pandharipande in \cite{pandharipande} and then by A. Schmitt in \cite{schmitt}. A recent calculation of their Brauer groups was carried out by R. Fringuelli and R. Pirisi in \cite{fringuelli}. However, much of their birational properties still remain a mystery.\\
We show that these moduli spaces are unirational for all ranks and degrees for genus $2\leq g\leq 9$ and for all ranks and fixed determinants for genus $g\geq 3$ for which $M_g$ is unirational which is known for $g \leq 14$. The proofs presented rely heavily on results about unirationality of $M_{g,g}, M_g$ and $P_{d,g}$ for different values of $g$. A summary of these results can be found in \cite{cas}, \cite{keneshlou}, \cite{verra} and \cite{bini}. \\
We shall extend Bertram's idea (in \cite{bertram}) of using a Poincar\'e bundle to show unirationality of the moduli space of semistable bundles over a fixed curve in rank $2$  to higher ranks in \S \ref{2.1s} and over families in \S\ref{bundle}. Lemma \ref{2.1} and the families we consider to prove the aforementioned results ensure that our proofs work for all degrees due to the presence of enough sections.\\
The families considered to show unirationality for universal moduli space of semistable bundles over smooth curves are different for genus $g = 2$ and $3\leq g\leq 9$. The primary reason behind this is the fact that generic curves of genus $g>2$ do not admit any non-trivial automorphisms. Hence, the open sublocus of automorphism free curves admits a fine moduli space with a universal family. In order to apply Lemma \ref{2.1}, we need $g$ general sections. This is achieved by considering the $g$-th fiber product of the universal curve over its parameter space. A corresponding such construction for the genus $2$ case is difficult to obtain directly from $M_2$ since the universal curve over it has general fiber $\p^1$ due to the presence of the hyperelliptic involution. The family considered for the genus $2$ case is thus constructed from a different approach- we utilise the fact that any genus $2$ curve with $2$ marked points can be realized as a plane curve of degree $4$ with exactly one node. We control the images of these marked points and these end up producing the two general sections as desired.\\
The parametrising families used here can be replaced with the ones described by Verra in \cite{verra} with similar proofs as discussed in \S\ref{4.3}.

\begin{ack}
\normalfont{I would like to thank Elham Izadi for suggesting this problem and for her continued guidance. I would also like to thank Alessandro Verra for pointing out the constructions described in \S\ref{4.3}.
}
\end{ack}

\section*{Notation}
\label{notation}
\begin{itemize}
\item We work over the field $\mathbb{C}$ of complex numbers
\item For an $A$-module $M$, $M^{\vee_A} = Hom(M,A)$. ($A$ being a commutative ring with unity)
\item For an $A$-module $M$, $\tilde{M}$ is the quasi-coherent sheaf on $Spec A$ associated to $M$
\item Curves are irreducible, smooth and projective varieties of dimension $1$ unless specified otherwise, family of curves refers to a smooth morphism of relative dimension $1$ with connected fibers.
\item $\p(W)$ is the space of $1$-dimensional subspaces of the vector space $W$
\item $Proj_S(\mathcal{E}) = Proj_S (Sym (\mathcal{E}))$ is the projective bundle of hyperplanes in $\mathcal{E}$, fibered over $S$
\item $M_g$ is the moduli space of smooth genus $g$ curves and $\mathcal{M}_g$ is the corresponding stack
\item $M_g^0$ is the moduli space of automorphism-free genus $g$ curves
\item $P_{d,g}$ is the universal Picard variety of degree $d$ line bundles over $M_g$ and $Pic_{d,g}$ is the universal Picard stack over $\mathcal{M}_g$
\item $U_C(r, {L})$ is the moduli space of semistable vector bundles of rank $r$ over $C$ with determinant ${L}$
\item $U_C(r,d)$ is the moduli space of semistable vector bundles of rank $r$ and degree $d$ over the curve $C$
    \item $U(r,d,g)$ is the universal moduli space of semistable vector bundles of rank $r$, degreee $d$ over smooth curves of genus $g$
    \item $U(r,\mathcal{L}, g)$ is the closed subscheme of $U(r,d,g)$, consisting of semistable vector bundles over smooth curves with determinant $\mathcal{L}_C\forall [C]\in M_g$, where $\mathcal{L}$ is a section of $Pic_{d,g}\rightarrow M_g$ 
\end{itemize}
Since $U(r,d,g)\cong U(r, d+r(2g-2), g), U(r,\mathcal{L}, g)\cong U(r, \mathcal{L}\otimes \omega^r, g)$, we shall assume $d> r(2g-1)$ and $deg(\mathcal{L})> r(2g-1)$ while proving unirationality of these moduli spaces. 

\section{Preliminaries}
\label{prelim}
We'll use the following lemma from \cite{verra}(Lemma 1.6):
\begin{lemma}
\label{2.1}
$C$ be a curve of genus $g$. Let $n_1,\cdots, n_g\in \Z$ be non-zero integers such that $n_1+\cdots+n_g = d$, then the map $a_{n_1,\cdots,n_g}: C^g\rightarrow Pic^d(C)$ given by $(x_1,\cdots, x_g)\mapsto n_1x_1+\cdots+n_gx_g$ is surjective.
\end{lemma}
\noindent We will need the following algebraic lemma, the proof of which is elementary and has been omitted here.
\begin{lemma}
\label{2.2}
    Let $M$ be a finite free $A$-module, $\mathfrak{m}$ a maximal ideal in $A$. Let $V = M\otimes A/\mathfrak{m}$, be the corresponding $(A/\mathfrak{m}=)$$k$-vector space. We have the following canonical isomorphisms:
    \begin{enumerate}
        \item $M^{\vee_A}\otimes_A A/\mathfrak{m}\cong V^{\vee}$ 
        \item Let $[id_M]\in Hom(M,M)\cong M\otimes_A M^{\vee_A}$. Then $[id_M]\mapsto [id_V]$ under the morphism $M\otimes_A M^{\vee_A}\xrightarrow{\otimes_A A/\mathfrak{m}} V\otimes_k V^\vee \cong Hom(V,V)$
    \end{enumerate}
\end{lemma}
\noindent Before we move on to Bertram's Poincar\'e bundle, it is crucial to understand how one can construct vector bundles in terms of extensions - using given bundles $\mathcal{E}_1,\mathcal{E}_2$; we construct a bundle $\mathcal{E}$ that fits into the exact sequence $0\rightarrow \mathcal{E}_1\rightarrow \mathcal{E}\rightarrow \mathcal{E}_2\rightarrow 0$. Any such bundle $\mathcal{E}$ is called an extension of $\mathcal{E}_2$ by $\mathcal{E}_1$.\\
    The equivalence classes of extensions of $\mathcal{E}_2$ by $\mathcal{E}_1$ (over a variety $X$) are given by elements of $H^1(X, Hom(\mathcal{E}_2, \mathcal{E}_1))\cong H^1(X, \mathcal{E}_2^\vee\otimes \mathcal{E}_1)$. The zero element of which corresponds to the trivial extension.\\
It is easy to see that non-zero scaling of vectors in $H^1(X, Hom(\mathcal{E}_2, \mathcal{E}_1))$ correspond to isomorphic extensions up to scaling of the sequence. Hence, isomorphism classes of (non-trivial)extension bundles could be realized as points in $\p H^1(X, Hom(\mathcal{E}_2, \mathcal{E}_1))$. \\
The above, coupled with the following from \cite{newstead}($137$, pg.$109$ and Lemma $5.2$, pg.$107$) provides enough background for Bertram's construction of a Poincar\'e bundle to conclude unirationality of certain moduli spaces.
\begin{lemma}
    \label{3.3}
    Let $F$ be a semistable bundle over a curve $C$ of rank $r$ and degree $d$, suppose $d> r(2g-1)$. Then we have:
    
    \begin{enumerate}
        \item $F$ is generated by its global sections
            
        \item $H^1(F) = 0$
    \end{enumerate}
 
    Furthermore, for a bundle $F$ over $C$ generated by its global sections, we have the following exact sequence: $$0\rightarrow \oplus^{r-1}\mathcal{O}_C\rightarrow F\rightarrow \det F\rightarrow 0$$
\end{lemma}
\subsection{Bertram's Poincar\'e Bundle}\hfill\\
\label{2.1s}
We give a quick presentation of the bundle considered by Bertram in \cite{bertram}(Definition-Claim $3.1$, pg.450).\\
We begin with a quick proof for unirationality of $U_C(2,L)$. \\
The key idea is to realize semistable bundles of high enough degree as extension classes of sequences dependent on the rank and degree.\\
Let $\p_L := \p(H^1(C,L^\vee))$, we construct a rational dominant map $\p_L \DashedArrow[->,densely dashed] U_C(2,\mathcal{L})$.\\
In order to do this, we construct a Poincar\'e Bundle on $C\times \p_L$ and consider the open semistable locus to define the map. We want an extension $$0\rightarrow \pi_{\p_L}^*\mathcal{O}_{\p_L}(1)\rightarrow \_\_ \rightarrow \pi_C^*L \rightarrow 0.$$
In order to do that, we consider the extension class 
\begin{align*}
[id]\in Hom(H^1(C, L^\vee), H^1(C, L^\vee))\cong H^1(C, L^\vee)^\vee\otimes H^1(C, L^\vee)\cong H^0(\p_L, \mathcal{O}_{\p_L}(1))\otimes H^1(C, L^\vee)\\
\underset{\text{Kunneth Formula}}{\cong}H^1(C\times \p_L,\pi_{\p_L}^*\mathcal{O}_{\p_L}(1)\otimes \pi_C^*L^\vee )
\end{align*}
This class corresponds to a Poincar\'e Bundle $\zeta$ on $C\times \p_L$ such that $\forall  p\in \p_L$, we have $[\zeta_p] = p$ as extension classes, where $\zeta_p :=  \zeta|_{C\times p}$.\\
The bundle $\zeta$ induces the desired rational map over the open locus of semistable bundles fibered over points in $\p_L$. The map is surjective by Lemma \ref{3.3}
\\
We can now discuss how Bertram's Poincar\'e Bundle could be adapted for higher ranks.\\
Analogous to the rank $2$ case, we consider a line bundle ${L}$ of degree $d$ on a smooth curve $C$ of genus $g$.\\
We are now ready to show that $U_C(r, L)$ is unirational $\forall r\geq 2$. \\
Let $\p_{r,L} = \p (\oplus^{r-1}H^1(C, {L^\vee}))$.\\
We would like to consider the following exact sequence on $C\times \p_{r,{L}}$: $$0\rightarrow \oplus^{r-1}\pi_{\p_{r,{L}}}^*\mathcal{O}_{\p_{r,L}}(1)\rightarrow \_\_ \rightarrow \pi_C^*{L}\rightarrow 0  \quad \quad \quad\quad (\star)$$
We would like to find a suitable extension class $[\mathcal{E}]$ for this sequence, such that for any $(p: \mathbb{C} \hookrightarrow\oplus^{r-1}H^1(C, {L^\vee}) )\in \p_{r,L}$, the restriction over $p$ $$0\rightarrow \oplus^{r-1}\mathcal{O}_C\rightarrow\mathcal{E}|_{C\times p}\rightarrow L \rightarrow 0 $$ corresponds to the extension class $[\mathcal{E}|_{C\times p}] = p(1) \in\oplus^{r-1}H^1(C, {L^\vee})$.\\
The space of extension classes is $$H^1(C\times \p_{r,L}, \oplus^{r-1}\pi_{\p_{r,L}}^*\mathcal{O}_{\p_{r,L}}(1)\otimes \pi_C^*{L}^\vee)\cong\oplus^{r-1}H^1(C,{L}^\vee)\otimes H^0(\p_{r,L}, \mathcal{O}_{\p_{r,L}}(1))= $$$$\oplus^{r-1}(H^1(C,{L}^\vee)\otimes (\oplus^{r-1}H^1(C,{L}^\vee)^\vee))\cong \oplus^{r-1}(\oplus^{r-1} (H^1(C,{L}^\vee)\otimes H^1(C,{L}^\vee)^\vee))\cong$$$$ \oplus^{r-1}(\oplus^{r-1}(End(H^1(C,L^\vee)))).$$ 
We consider the class $\zeta\in \oplus^{r-1}(\oplus^{r-1}(End(H^1(C, L^\vee))))$ where \begin{equation}\label{zeta}\zeta = \left[\begin{matrix}
    [id] &0 &0 &\cdots &0\\ 0& [id]&0&\cdots &0\\   \vdots & \dots&\dots&\ddots & \vdots\\
    0 & 0&  \dots  & 0&[id]
\end{matrix}\right].\end{equation}
Let $\mathcal{E}$ be the corresponding bundle. We will show that $\mathcal{E}$ satisfies the property described above. Let $p:\mathbb{C}\hookrightarrow{}\oplus^{r-1}H^1(C,L^\vee)$ be a point in $\p_{r,L}$ with $p(1) = (v_1,\cdots, v_{r-1})$. We have
$$[\mathcal{E}|_{C\times p}] \in H^1(C\times p, \oplus^{r-1}Im(p)^\vee\otimes {L}^\vee) \cong \oplus^{r-1}H^1(C,{L}^\vee) .$$ Let $\zeta|_{C\times p}  = ({}^{p}\zeta_i)_i$. So, ${}^p\zeta_i = \zeta_i.p(1) = v_i$ where $\zeta_i$ is the $i$-th column of the matrix $\zeta$. Thus, $[\zeta|_{C\times p}] = (v_i)_i = p(1)$ as desired.
This induces a rational surjective morphism \begin{equation}\label{phi}
    \phi: \p_{r,L}\DashedArrow[->,densely dashed] U_C(r,{L})
\end{equation} defined over the open locus of $\p_{r,L}$ parametrising semistable bundles over $C$ by Lemma \ref{3.3}.
\subsection{Bertram's extension class over families}\hfill\\
\label{bundle}
The expression for the extension class $[\zeta]$ in (\ref{zeta}) allows for a similar choice of an extension class $[\delta]$ over certain families of curves so that the restriction of $\delta$ to each curve in this family is equal to the extension class $[\zeta]$. \\
Let $\mathscr{C}\rightarrow U = Spec(A)$ be a family of curves over $U$ of genus $g$ and $\mathcal{L}\rightarrow\mathscr{C}$ be a line bundle over $\mathscr{C}$ of relative degree $d$. Fix an $r\geq 2$ (as explained in \S\ref{notation}, we assume $d>r(2g-1)$).\\
By Cohomology and Base Change, we have that $R^1\pi_*\mathcal{L}^\vee$ is locally free and there are canonical isomorphisms $(R^1\pi_*\mathcal{L}^\vee)_u \cong H^1(\mathscr{C}_u, \mathcal{L}|_{\mathscr{C}_u})\forall u\in U$. Since $\mathscr{C}\rightarrow U$ is flat and $\forall u\in U$, we have $h^1(C,L^\vee) = h^0(C, L+K_C) = d+ g-1 $ where $\mathscr{C}_u = C, \mathcal{L}|_{\mathscr{C}_u} = L$. \\ Let $\mathcal{E} :=\oplus^{r-1}R^1\pi_*\mathcal{L}^\vee $.
We show the existence of a Poincar\'e Bundle $\zeta$ on $\p\times_{U} \mathscr{C}$ so that $\zeta$ restricts to Bertram's Poincar\'e Bundle over all points $u\in U$ as in \S \ref{2.1s} where $\p = \mathbb{P}_U (\mathcal{E})$.
We have the following diagram:
$$\begin{tikzcd}
&\zeta\arrow{d}\\
&\p\times_U \mathscr{C}\arrow{d}\\
\mathcal{L}\arrow{d}& \p \arrow{d}{\pi_U}\\
\mathscr{C}\arrow{r}{\pi}&U\end{tikzcd}
$$\\
In order to construct $\zeta$, we consider the extension classes for the sequence: $$0\rightarrow \oplus^{r-1}\pi_{\p}^* \mathcal{O}_{\p}(1)\rightarrow \_\_\_\rightarrow\pi_\mathscr{C}^* \mathcal{L}\rightarrow 0$$
The extension space is: $H^1(\p\times_{U}\mathscr{C}, \oplus^{r-1}\pi_{\p}^* \mathcal{O}_{\p}(1)\otimes \pi_\mathscr{C}^* \mathcal{L}^\vee)$.
We consider the following class: \\
$$[\zeta] := \left[\begin{matrix}
    [id] &0 &0 &\cdots &0\\ 0& [id]&0&\cdots &0\\   \vdots & \dots&\dots&\ddots & \vdots\\
    0 & 0&  \dots  & 0&[id]
\end{matrix}\right]\in \oplus^{r-1}(\oplus^{r-1}(End(H^1(\mathscr{C},\mathcal{L^\vee}))) \cong \oplus^{r-1}((H^1(\mathscr{C},\mathcal{L}^\vee)^{\vee_A})^{\oplus r-1}\otimes H^1(\mathscr{C},\mathcal{L}^\vee))$$
Shrinking $U$ if necessary, let's additionally assume that $\mathcal{E}$ is free. We have by \cite{har} (Prop \RomanNumeralCaps{3}, 8.5), $\mathcal{E} \cong \oplus^{r-1} (H^1(\mathscr{C}, \mathcal{L^\vee}))^\sim \implies \mathcal{E^\vee} \cong  (\oplus^{r-1}H^1(\mathscr{C}, \mathcal{L^\vee})^{\vee_A})^\sim$ and  $ H^0(\p, \mathcal{O}_{\p}(1))\cong H^0(U, \pi_{U *}\mathcal{O}_{\p}(1)) = H^0(U, \mathcal{E^\vee})\cong \oplus^{r-1} H^1(\mathscr{C}, \mathcal{L^\vee})^{\vee_A}$.\\
Thus, we have the class: $$[\zeta]\in \oplus^{r-1}(H^0(\p, \mathcal{O}_{\p}(1)))\otimes H^1(\mathscr{C},\mathcal{L}^\vee))\xrightarrow{\cup}\oplus^{r-1} H^1(\p\times_{U}\mathscr{C},\pi_{\p}^* \mathcal{O}_{\p}(1)\otimes \pi_\mathscr{C}^* \mathcal{L}^\vee ) $$$$\cong H^1(\p\times_{U}\mathscr{C}, \oplus^{r-1}\pi_{\p}^* \mathcal{O}_{\p}(1)\otimes \pi_\mathscr{C}^* \mathcal{L}^\vee)$$
We show that this is the desired extension class.
For any $u\in U$, by cohomology and base change, we have $$\mathcal{E}_u^\vee = \oplus^{r-1}(R^1\pi_{u*}L^\vee)^\vee\cong \oplus^{r-1}H^1(C, L^\vee)^\vee\implies \p_u = \p_{r,L}.$$ 
Let $ C = \mathscr{C}_u$ and $L = \mathcal{L}|_{\mathscr{C}_u}$, thus $(\p'\times_{U'}\mathscr{C}')_u = \p_L\times C$. We have the following diagram:$$
\begin{tikzcd}
    L \arrow{d}&\p_u = \p_{r,L
    }\arrow[d]\\
    C \arrow{r}{\pi_u}& u
\end{tikzcd} $$
\noindent The restriction $\zeta_u$ corresponds to an extension class of the following sequence:
$$0\rightarrow \oplus^{r-1}\pi_{\p_u}^*\mathcal{O}_{\p_u}(1)\rightarrow \_\_\rightarrow \pi_C^*L\rightarrow 0
$$ which is the image of $[\zeta]$ under the map: $$\oplus^{r-1}(\oplus^{r-1}(End(H^1(\mathscr{C}',\mathcal{L^\vee}|_{\mathscr{C}'}))) \rightarrow H^1(\p'\times_{U'}\mathscr{C}', \oplus^{r-1}\pi_{\p'}^* \mathcal{O}_{\p'}(1)\otimes \pi_\mathscr{C'}^* \mathcal{L}^\vee) \rightarrow $$$$H^1((\p'\times_{U'}\mathscr{C}')_u, \oplus^{r-1}(\pi_{\p'}^* \mathcal{O}_{\p'}(1)\otimes \pi_\mathscr{C'}^* \mathcal{L}^\vee)_u) = H^1(C\times \p_{r,L}, \oplus^{r-1}\pi_{\p_{r,L}}^*\mathcal{O}_{\p_{r,L}}(1)\otimes \pi_C^*{L}^\vee) $$$$= \oplus^{r-1}(\oplus^{r-1}(End(H^1(C,{L}^\vee)))$$\\
Now, by Lemma \ref{2.2}, the identity endomorphism on $H^1(\mathscr{C}', \mathcal{L}^\vee|_{\mathscr{C}'})$ restricts to the identity endomorphism on $H^1(C, L^\vee)$ and we have that $[\zeta_u]$ is the same extension class as the one used in \S \ref{2.1s} to construct Bertram's Poincar\'e Bundle as claimed.
\section{Unirationality of $U(\lowercase{r,d,g})$}
\label{4s}
\subsection{Genus $2$}\hfill\\
\label{3.1}
We look at genus $2$ curves and their realisations as nodal plane curves of degree $4$.\\
Let $C$ be a curve of genus $2$ and $p_1,p_2\in C$ such that $3p_1+p_2\not\sim 2K_C$, which is satisfied by a general tuple $(C,p_1,p_2)\in M_{2,2}$.\\
The image of the map $f_{p_1,p_2}: C \rightarrow \p^2$ given by $|3p_1+p_2|$ is then a plane curve of degree $4$ with exactly one node.
Blowup of this image at the node recovers the curve $C$.\\
We fix distinct and non-colinear points $P_1, P_2, P_3\in \p^2$. We consider the space $S (\cong\p^7)$ of degree $4$ plane curves that pass through $P_1, P_2, P_3$ with a node at $P_3$ and a $3$-flex at $P_1$ which passes through $P_2$.\\ 
Let $U\underset{\text{open}}{\subset} S$ be the collection of nodal curves in $S$ that are smooth away from $P_3$.\\
Any $u\in U$ can be realized as the image of a genus $2$ curve $C$ under the map $f_{p_1,p_2}:C \rightarrow\p^2$ for some $x_1,x_2\in C$ with $f(p_1) = P_1, f(p_2) = P_2$.\\
Let $\mathcal{N}\xrightarrow{\pi} U$ be the universal family of these nodal curves over $U$. We have $U\times P_i \subset \mathcal{N}\subset U\times \p^2 \forall i\in \{1,2,3\}$.\\
Let $s_i: U\rightarrow \mathcal{N}$ be the section given by $P_i$. Let $\mathscr{C} = Bl_{s_3}\mathcal{N}$ be the blowup of $\mathcal{N}$ along the section $s_3$.\\
Along the fibers of $\pi$, this blowup is the blowup of a nodal curve at its node - which results in its desingularization as described in \S \ref{prelim}. Thus, we get $\mathscr{C}\rightarrow\mathcal{N}\xrightarrow{\pi}U$ is a family of genus $2$ curves.\\
Since the sections $s_1,s_2$ do not intersect $s_3$, they lift to $\mathscr{C}$. \\
So we have $
\psi:\mathscr{C}\rightarrow U $ with sections $s_1, s_2: U\rightarrow \mathscr{C}$.
A general tuple $(C,p_1,p_2)\in M_{2,2}$, can be realized as $(\psi^{-1}(u), s_1(u), s_2(u))$ for some $u\in U$ since Im$(f_{x_1,x_2})\subset \p^2$ has a linear automorphism which sends the node to $P_3$ and $x_i$ to $P_i$ ($i= 1,2$).\\
For any pair of non-zero integers $n_1,n_2$ such that $n_1+n_2 = d (> 3r)$, the line bundle $\mathcal{L}= \mathcal{O}_{\mathscr{C}}(n_1s_1+n_2s_2)$ induces a morphism $U \rightarrow P_{d,2}$. By \S\ref{2.1s}, the image of this morphism is dense;  giving us the following diagram:$$
\begin{tikzcd}
    \mathcal{L} = \mathcal{O}_\mathscr{C}(n_1s_1+n_2s_2)\arrow{d}\\
    \mathscr{C}\arrow{r}& U\arrow[r]
    &P_{d,2}
\end{tikzcd}$$
We're now ready to prove unirationality for $U(r,d,2)$. \\
Following the notation of \S\ref{bundle}, we shall prove unirationality by constructing a rational map $\p\DashedArrow[->,densely dashed] U(r,d,2)$ with a dense image. In order to do this, we shall move our attention to an affine open $U'=\mathrm{Spec}(A)\subset U$ so that $\mathcal{E}|_{U'}$ is free and 
the restriction $U'\rightarrow P_{d,2}$ is still dominant. By \S\ref{bundle}, we have a Poincar\'e Bundle $\zeta$ on $\p'\times_{U'} \mathscr{C}'$ so that $\zeta$ restricts to Bertram's Poincar\'e Bundle over all points $u\in U'$ as presented in \S \ref{2.1s}.
We consider the following maps:\begin{equation}
\begin{tikzcd}
\label{diagram:main}&\zeta\arrow{d}\\
&\p'\times_{U'}\mathscr{C}'\arrow{d}\\
    \mathcal{L} = \mathcal{O}_\mathscr{C}(n_1s_1+n_2s_2)\arrow{d} & \p'\subset\p \arrow{d}{\pi_U}\arrow[r, dashed]{}{\zeta_{ss}} & U(r,d,2)\\
    \mathscr{C}\arrow{r}{\pi}& U
\end{tikzcd}\end{equation}
The rational morphism induced by $\zeta$, denoted by
$\Phi:\p\DashedArrow[->,densely dashed]  U(r,d,2)$, is defined over the open locus of points in $\p'$ over which $\zeta$ restricts to a semistable bundle.\\
Thus, we have the following diagram:
$$\begin{tikzcd}
    \p'\arrow{d}\arrow[r, "\zeta_{ss}", dashed]&U(r,d,2)\arrow{d}{\det}\\
    U'\arrow{r}{\theta}& P_{d,2}
\end{tikzcd}$$The rational map $\zeta_{ss}$ is defined on the points $p\in \p'$ for which $\zeta|_{C_p\times p}$ is semistable.\\
Now for any point $q = (L\rightarrow C)\in P_{d,2}$, the fiber over $q$ under $det$ is $ U_C(r, L)$. If $q=\theta(u)$ for some $u\in U'$, then we have that the restriction $\p'_u\DashedArrow[->,densely dashed] U_C(r,L)$ is the rational map given by Bertram's Poincar\'e Bundle as shown above and thus, is surjective.\\
Since $\theta$ is a dominant morphism, so is the rational map $\p'\DashedArrow[->,densely dashed] U(r,d,2)$ and since $\p' \underset{\text{open}}{\subset}\p$, we have the desired dominant rational map $\p\DashedArrow[->,densely dashed]  U(r,d,2)$ proving unirationality.

\subsection{Genus $3\leq \lowercase{g}\leq 9$}\hfill\\
\label{3.2}
If $g\geq 3$, a general curve $C\in M_g$ is known to not have any non-trivial automorphisms and the open sub-locus of these automorphism free curves is denoted by $M_g^0$ which is a fine moduli space.\\
Let $\mathscr{C}_g^0$ be the universal family over this open sub-locus $M_g^0$.\\
Set $\mathscr{C}_{g}^{0,n} := \times_{M_g^0}^n \mathscr{C}_g^0 \forall n \geq 1$. We have $$\begin{tikzcd}
 \mathscr{C}_{g}^{0,n+1} \arrow{d}{\psi_{n+1}}&(C, p_1,\cdots, p_{n+1})\arrow[mapsto]{d}\\\mathscr{C}_{g}^{0,n}& (C, p_1,\cdots,p_n)   
\end{tikzcd}$$ is the pullback of the universal curve to $\mathscr{C}_{g}^{0,n}$.
It is easy to see that $\mathscr{C}_{g}^{0,g}$ and $ M_{g,g}$ are birational.
We will show that there is a dominant morphism from $\mathscr{C}_{g}^{0,g}$ to $P_{d,g}\forall d$.\\
We have $g$ sections $\mathscr{C}_{g}^{0,g}\xrightarrow{\sigma_i}\mathscr{C}_{g}^{0,g+1}$ given by $(C, p_1,\cdots, p_g)\mapsto (C, p_1,\cdots, p_g, p_i)\forall 1\leq i\leq g$. Since $d> r(2g-1) > g, \exists n_1,\cdots, n_g \in \mathbb{Z}_{>0}$ such that $n_1+\cdots+n_g = d$. \\Define the line bundle $\mathcal{L} = \mathcal{O}_{\mathscr{C}_{g}^{0,g+1}}(\Sigma n_i \sigma_i)$ of degree $d$ $$
\begin{tikzcd}
    \mathcal{L} = \mathcal{O}_{\mathscr{C}_{g}^{0,g+1}}(\Sigma n_i \sigma_i)\arrow{d}\\
    \mathscr{C}_{g}^{0,g+1}\arrow{r}{\psi_{n+1}}&
    \mathscr{C}_{g}^{0,g}
\end{tikzcd} $$ 
By Lemma ~\ref{2.1}, we have that the image of the induced morphism $\mathscr{C}_{g}^{0,g}\xrightarrow{\pi_g} P_{d,g}$ contains all the points $\left\{
    L\rightarrow C
|\forall [C]\in M_g^0, deg(L) = d\right\} = f^{-1}(M_g^0)$ where $f:P_{d,g}\rightarrow M_g$ is the forgetful functor, showing that $\pi_g$ is dominant.\\
From \cite{keneshlou}(\S 3), we have $M_{g,g}$ is unirational $\forall g$ with $3\leq g \leq 9$, thus we have that $\mathscr{C}_{g}^{0,g}$ is unirational as well.\\
We are now ready to show unirationality of $U(r,d,g)$.\\
Since $\mathscr{C}_{g}^{0,g}$ is unirational, $\exists \phi:\p^N \DashedArrow[->,densely dashed] \mathscr{C}_{g}^{0,g}$ with a dense image, let $\phi$ be defined over $U\underset{\text{open}}{\subset} \p^N$. Since $\pi_g$ is dominant, we have that the composition $\Phi: U\rightarrow P_{d,g}$ is dominant as well.\\
Let $\mathcal{E} := \oplus^{r-1}R^1\pi_{U,*}\mathcal{L}^\vee$. Shrinking $U$ if necessary, we can additionally assume that $U$ is affine and $\mathcal{E}$ is free. Following the same notation as in \S\ref{bundle}, the pullback of $(\mathcal{L}\rightarrow \mathscr{C}_g^{0,g+1})$ over $U$ is given by the following diagram:
$$
\begin{tikzcd}
    \mathcal{L}\arrow{d}& \mathcal{L}_U\arrow{d}& \p = \mathbb{P}_U(\mathcal{E})\arrow{d}\\ 
    \mathscr{C}_{g}^{0,g+1}\arrow{rd}& \arrow{l}\mathscr{C}_U \arrow{r}{\pi_U}& U\arrow{ld}\\
    & \mathscr{C}_{g}^{0.g}
\end{tikzcd}$$
By \S\ref{bundle}, we have a Poincar\'e Bundle $\zeta$ on $\p\times_{U}\mathscr{C}$ so that
restriction of this bundle over any point $u\in U$ is Bertram's Poincar\'e Bundle as considered in \S \ref{2.1s}.\\ We thus have the following diagram as in (\ref{diagram:main}):

$$\begin{tikzcd}
\zeta\arrow{d}&\\
\p\times_{U}\mathscr{C}\arrow{d}&\\
  \p\arrow{d}\arrow[r, "\zeta_{ss}", dashed]&U(r,d,g)\arrow{d}{\det}\\
    U\arrow{r}{\Phi}& P_{d,g}
\end{tikzcd}$$The rational map $\zeta_{ss}$ is defined over the points $p\in \p$ for which $\zeta|_{C_p\times p}$ is semistable where $C_p = \pi_U^{-1}(p)$.\\
Now for any point $q = (
    L, C
)\in P_{d,g}$, the fiber over $q$ under $det$ is $U_C(r, L)$. If $q=\Phi(u)$ for some $u\in U$, then we have that the restriction $\p_u\DashedArrow[->,densely dashed] U_C(r,L)$ is canonically isomorphic to the rational map in (\ref{phi}) and thus, is surjective.\\
Since $\Phi$ is dominant, we have $\p\DashedArrow[->,densely dashed] U(r,d,g)$ is dominant as well. Now, $\p$ is a projective bundle over $U\underset{\text{open}}{\subset}\p^N$, we have that $\p$ is a rational variety. Thus, $U(r,d,g)$ is unirational $\forall r,d,g$ with $2\leq r, 3\leq g\leq 9$.
\begin{rem} {\normalfont
We know that $\kappa(P_{d,g})\geq 0$ for $g\geq 10$ from \cite{martin} and \cite{bini}. Therefore, $U(r,d,g)$ is unirational if and only if $2\leq g \leq 9$.
}\end{rem}
\subsection{Descriptions by Mukai and Verra}\hfill\\
\label{4.3}
Verra, in \cite{verra}, used certain rational homogeneous spaces constructed by Mukai in \cite{mukai} to construct certain varieties $P_g\subset \p^{\dim P_g + g-2}$ with the property that a general curve of genus $g$ can be realized as a curvilinear section of $P_g$. \\
These $P_g$ can then be used to show unirationality of $M_{g,g}$ for $7\leq g\leq 9$ which further implied unirationality for $P_{d,g}$. The case for $4\leq g \leq 6$ works a little differently. \\
We still use $P_g$ with the same property that a general curve of genus $g$ can be realized as a curvilinear section of it but we don't have $P_g\subset \p^{\dim P_g + g-2}$ anymore:\\
For genus $4$, $P_4$ is a general complete intersection of type $(2,3)$ in $\p^6$.\\
For genus $5$, $P_5$ is a fourfold (general) complete intersection of type $(2,2,2)$ in $\p^7$.\\
For genus $6$, $P_6$ is a fivefold which is a general quadratic section of $Gr(2,5)$.
\section{Unirationality of $U(\lowercase{r},\mathcal{L} , \lowercase{g})$}
\label{4}
We shall prove that for any section $\mathcal{L}: \mathcal{M}_g\rightarrow Pic_{d,g}$ of the forgetful morphism $Pic_{d,g}\rightarrow M_g$, unirationality of $U(r,\mathcal{L}, g)$ and $\mathcal{M}_g$ are equivalent. A summary of recent developments on unirationality of $\mathcal{M}_g$ can be found in \cite{farkas} and \cite{keneshlou}.\\
Some examples are: $\mathcal{O}: M_g\rightarrow Pic_{0,g}$, given by $[C]\mapsto 
    (\mathcal{O}_C\rightarrow C)$ and $\omega^{\otimes n}:M_g\rightarrow Pic_{n(2g-2), g}$, given by $[C]\mapsto (
    K_C^{\otimes n}\rightarrow C)
, n\in \mathbb{Z}$. In fact, these are the only possible sections by the strong Franchetta Conjecture \cite{mestrano}.\\
Now suppose $\mathcal{M}_g$ is unirational for some $g$ and we have a section $\mathcal{L}:\mathcal{M}_g\rightarrow Pic_{d,g}$. Thus, we have a rational map $\Phi: \p^N\DashedArrow[->,densely dashed]\mathcal{M}_g$ for some $N\in \mathbb{N}$ with dense image.\\
Hence, $\exists U\underset{\text{open}}{\subset}\p^N$ such that $\Phi: U\rightarrow \mathcal{M}_g$ is a dominant morphism.
The morphism $\Phi$ induces a curve of genus $g$, $\mathscr{C}\xrightarrow{\pi} U$. $\mathcal{L}\circ \Phi$ induces a morphism $U\rightarrow Pic_{d,g}$ which also induces a line bundle $(
    \mathcal{L}\rightarrow\mathscr{C})$ of degree $d$ over this curve.
Let $\mathcal{E}:=\oplus^{r-1}R^1\pi_*\mathcal{L}^\vee $. Following the same notation as in \S\ref{bundle}, we have the following diagram:
$$\begin{tikzcd}
    \mathcal{L}\arrow{d}& \p = \mathbb{P}_U(\mathcal{E})\arrow{d}\\ \mathscr{C}\arrow{r}& U
\end{tikzcd}$$
By \S\ref{bundle}, we have a Poincar\'e Bundle on $\p'\times_{U'}\mathscr{C}'$ for some affine open $U'\subset U$ so that $U'\xrightarrow{\Phi}M_g$ is dominant and $\mathcal{E}|_{U'}$ is trivial where $\p' = \p|_{U'}$, $\mathscr{C}'=\mathscr{C}|_{U'}$.
Hence we have the following diagram:
$$
\begin{tikzcd}
    \p'\arrow{d}\arrow[r, "\zeta_{ss}", dashed]&U(r,\mathcal{L},g)\arrow{d}\\
    U'\arrow{r}{\Phi}& M_g
\end{tikzcd}$$The rational map $\zeta_{ss}$ is defined on the points $p\in \p'$ for which $\zeta|_{C_p\times p}$ is semistable.\\
Now, for any point $q = [C]\in M_g$, the fiber in $U(r, \mathcal{L},g)$ is given by $ U_C(r, L)$. Restriction of $\p'$ over any point $u\in U'$ is canonically isomorphic to the projective bundle considered in \S \ref{2.1s}. If $q=\Phi(u)$ for some $u\in U'$, then we have that the restriction $\p'_u\DashedArrow[->,densely dashed] U_C(r,L)$ is canonically isomorphic to the map (\ref{phi}) and thus, is surjective.\\
Since $\Phi$ is a dominant, we can say the same about the rational map $\p'\DashedArrow[->,densely dashed] U(r,\mathcal{L},g)$. Now, $\p$ is a projective bundle over $U\underset{\text{open}}{\subset}\p^N$, we have that $\p$ is a rational variety. Thus, $U(r,\mathcal{L},g)$ is unirational.

\section*{references}

\renewcommand{\section}[2]{}%

\bibliographystyle{amsalpha}
\bibliography{ack}

\providecommand{\bysame}{\leavevmode\hbox to3em{\hrulefill}\thinspace}
\providecommand{\MR}{\relax\ifhmode\unskip\space\fi MR }
\providecommand{\MRhref}[2]{%
  \href{http://www.ams.org/mathscinet-getitem?mr=#1}{#2}
}
\providecommand{\href}[2]{#2}
\begin{thebibliography}{GBV12}

\bibitem[AC81]{ac}
E.~Arbarello and M.~Cornalba, \emph{Footnotes to a paper of {B}eniamino {S}egre}, Mathematische Annalen 256 (1981), 341 -- 362.

\bibitem[Ber92]{bertram}
A.~Bertram, \emph{Moduli of rank-2 vector bundles, theta divisors, and the geometry of curves in projective space}, J. Differential Geom. \textbf{35} (1992), no.~2, 429 -- 469.

\bibitem[CF07]{cas}
G.~Casnati and C.~Fontanari, \emph{On the rationality of moduli spaces of pointed curves}, Journal of the London Mathematical Society \textbf{75} (2007), 582–596.

\bibitem[CR84]{cr1}
M.~C. Chang and Z.~Ran, \emph{Unirationality of the moduli space of curves of genus 11, 13 (and 12)}, Inventiones Math. \textbf{76} (1984), 41 --– 54.

\bibitem[EH87]{EH}
D.~Eisenbud and J.~Harris, \emph{The {K}odaira dimension of the moduli space of curves of genus $\geq 23$}, Inventiones Math. \textbf{90} (1987), 359–387.

\bibitem[FP]{fringuelli}
R.~Fringuelli and R.~Pirisi, \emph{The {B}rauer {G}roup of the {U}niversal {M}oduli {S}pace of {V}ector {B}undles {O}ver {S}mooth {C}urves}, IMRN \textbf{2021}, no.~18, 13609–13644.

\bibitem[FV]{farkas}
G.~Farkas and A.~Verra, \emph{On the {K}odaira dimension of the moduli space of curves of genus 16}.

\bibitem[GBV12]{bini}
C.~Fontanari G.~Bini and F.~Viviani, \emph{On the {B}irational {G}eometry of the {U}niversal {P}icard {V}ariety}, IMRN \textbf{2012} (2012), no.~4, 740--780.

\bibitem[GFP23]{FJP}
D.~Jensen G.~Farkas and S.~Payne, \emph{The kodaira dimension of {$\overline{\mathcal{M}}_{22}$} and {$\overline{\mathcal{M}}_{23}$}}.

\bibitem[Har77]{har}
R.~Hartshorne, \emph{Algebraic geometry}, Springer-Verlag, 1977.

\bibitem[KS99]{SK}
A.~King and A.~Schofield, \emph{Rationality of moduli of vector bundles on curves}, Indag. Math. \textbf{10} (1999), 519–535.

\bibitem[KT21]{keneshlou}
H.~Keneshlou and F.~Tanturri, \emph{On the unirationality of moduli spaces of pointed curves}, Math. Z. \textbf{299} (2021), 1973–1986.

\bibitem[KV17]{martin}
Sebastian Casalaina-Martin{,} Jesse~Leo Kass and Filippo Viviani, \emph{The singularities and birational geometry of the universal compactified jacobian}, Algebraic Geometry \textbf{4} (2017), no.~3, 353--393.

\bibitem[Mes87]{mestrano}
N.~Mestrano, \emph{Conjecture de {F}ranchetta forte}, Inventiones Mathematicae (1987), 365--376.

\bibitem[Muk88]{mukai}
S.~Mukai, \emph{Curves, {K}3 surfaces and {F}ano 3-folds of {G}enus {$\leq 10$}}, Algebraic Geometry and Commutative Algebra In Honor of Masayoshi Nagata \textbf{1} (1988), 357--377.

\bibitem[New75]{N}
P.E. Newstead, \emph{Rationality of moduli spaces of stable bundles}, Math. Ann. \textbf{215} (1975), 251–268.

\bibitem[New11]{newstead}
\bysame, \emph{Introduction to {M}oduli {P}roblems and {O}rbit {S}paces}, vol.~17, 2011.

\bibitem[Pan96]{pandharipande}
R.~Pandharipande, \emph{A {C}ompactification {O}ver {$\overline{M_g}$} {O}f {T}he {U}niversal {M}oduli {S}pace of {S}lope-{S}emistable {V}ector {B}undles}, Journal of the American Mathematical Society \textbf{9} (1996), no.~2, 425--471.

\bibitem[Sch04]{schmitt}
A.~Schmitt, \emph{The {H}ilbert compactification of the universal moduli space of semistable vector bundles over smooth curves}, J. Differential Geometry \textbf{66} (2004), no.~2, 169--209.

\bibitem[Ser81]{ser}
E.~Sernesi, \emph{L{'}unirazionalità della varietà dei moduli delle curve di genere 12}, Annali della Scuola Normale Superiore di Pisa \textbf{8} (1981), 405–439.

\bibitem[Sev15]{sev}
F.~Severi, \emph{Sulla classificazione delle curve algebriche e sul teorema d’esistenza di riemann}, Rendiconti della Reale Accademia Naz. Lincei \textbf{24} (1915), 877–888.

\bibitem[Ver]{verra}
A.~Verra, \emph{The {U}nirationality of the {M}oduli {S}pace of {C}urves of {G}enus 14 or {L}ower}, Compositio Mathematica \textbf{141}, no.~6, 1425–1444.

\end{thebibliography}

\end{document}